\documentclass[reqno,12pt]{amsart}

\def\3{\ss}
\def\({\left(}
\def\){\right)}
\def\[{\left[}
\def\]{\right]}

\def\po#1#2{(#1)_#2}
\def\poq#1#2{(#1;q)_#2}

\numberwithin{equation}{section}

\theoremstyle{remark}

\begin{document}

\newbox\Adr
\setbox\Adr\vbox{
\centerline{\sc C.~Krattenthaler$^\dagger$ and H.~Rosengren}
\vskip18pt
\centerline{Institut Girard Desargues,
Universit\'e Claude Bernard Lyon-I}
\centerline{21, avenue Claude Bernard, F-69622 Villeurbanne Cedex, France}
\centerline{e-mail: \tt kratt@euler.univ-lyon1.fr}
\centerline{WWW: \footnotesize\tt http://euler.univ-lyon1.fr/home/kratt}
\vskip10pt
\centerline{Department of Mathematics, Chalmers University of Technology}
\centerline{and G\"oteborg University, SE-412\,96 G\"oteborg, Sweden}
\centerline{e-mail: \tt hjalmar@math.chalmers.se}
\centerline{WWW: \footnotesize\tt http://www.math.chalmers.se/\~{}hjalmar}
}

\title{On a hypergeometric identity of Gelfand, Graev and Retakh}

\author[C. Krattenthaler and H. Rosengren]{\box\Adr}

\address{Institut Girard Desargues, Universit\'e Claude Bernard Lyon-I,
21, avenue Claude Bernard, F-69622 Villeurbanne Cedex, France.
E-mail: {\tt kratt@euler.univ-lyon1.fr},\newline
WWW: \tt http://euler.univ-lyon1.fr/home/kratt.}

\address{Department of Mathematics, Chalmers University of Technology
and G\"oteborg University, SE-412\,96 G\"oteborg, Sweden.
E-mail: {\tt hjalmar@math.chalmers.se},\newline
WWW: \tt http://www.math.chalmers.se/\~{}hjalmar.}

\thanks{$^\dagger$ Research partially supported by the Austrian
Science Foundation FWF, grant P12094-MAT, and by
EC's IHRP Programme, grant HPRN-CT-2001-00272.}



\begin{abstract}
A hypergeometric identity equating a triple sum to a single sum, originally
found by Gelfand, Graev and Retakh [Russian Math.\ Surveys 
{\bf 47} (1992), 1--88] by using systems of differential equations, is
given hypergeometric proofs. As a bonus, several $q$-analogues can be derived.
\end{abstract}

\maketitle

\begin{section}{Introduction}

In a recent talk at the Institute of Mathematics and Its Applications,
Minnesota, Richard Askey posed several problems for the audience, one
of which concerned the following identity due to
Gelfand, Graev and Retakh \cite[p.~67]{GeGRAA},
\begin{align}\label{GGR0}
\sum _{j,k,m\ge0} ^{}
&{\frac{(\alpha)_j\,(\beta)_k\,(1-\gamma)_m\,(\gamma)_{j+k-m}} 
{j!\,k!\,m!\,(\alpha+\beta)_{j+k-m}}}x^jy^kz^m\\
\label{GGR1}
&\kern1cm=\frac {(1-z)^{\alpha+\beta-1}\,
(1-xz)^{\gamma-\alpha-\beta}}
{(1-x)^{\gamma-\beta}\,(1-y)^{\beta}}
 {}_2F_1\left[\begin{matrix}{ \beta,\alpha+\beta-\gamma}\\{ \alpha+\beta
}\end{matrix};{\frac{(x-y)(1-z)} {(1-y)(1-xz)}}\right]\\
\label{GGR2}
&\kern1cm= \frac{{\left( 1 - z \right) }^{ \alpha + \beta-1} \,
     {\left( 1 - x z \right) }^{\gamma-\alpha } }{
     {\left( 1 - x \right) }^{\gamma} \,
     {\left( 1 - y z \right) }^{\beta}}
     {} _{2} F _{1} \!\left [ \begin{matrix} { \beta,
      \gamma}\\ { \alpha + \beta}\end{matrix} ;
      {\displaystyle \frac{\left( y - x \right)  \left( 1 - z \right)
      }{\left( 1 - x \right)  \left( 1 - y z \right) }}\right ] ,
\end{align}
where the Pochhammer symbol $(a)_k$ is defined by
$(a)_k  =  a(a+1)(a+2)\cdots (a+k-1)$ if
$k>0$, $(a)_0  = 1$, $(a)_{k}=1/(a+k)_{-k}$ if $k<0$,
and the hypergeometric series $_rF_s$
is defined by
\begin{equation*}  
{}_r F_s\!\left[\begin{matrix} a_1,\dots,a_r\\ b_1,\dots,b_s\end{matrix}; 
z\right]=\sum _{\ell=0} ^{\infty}\frac {\po{a_1}{\ell}\cdots\po{a_r}{\ell}}
{\ell!\,\po{b_1}{\ell}\cdots\po{b_s}{\ell}} z^\ell\ .
\end{equation*}
This triple sum -- single sum identity is derived in \cite{GeGRAA} as
a special case of a general reduction formula for hypergeometric functions
connected with the Gra\3mannian $G_{N,K}$, and is derived there by
exploiting systems of differential equations. Richard Askey posed
the problem of finding a purely hypergeometric proof of this identity
(i.e., one that uses classical summation and transformation formulas for
hypergeometric series), with the ulterior motive that such a proof will
make it possible to find a $q$-analogue of the identity, which was not
found until then.

The purpose of this note is to provide such a hypergeometric proof 
for this identity, which we give in the following section. 
Indeed, as a bonus we are able to 
derive several $q$-analogues of the formula, see Section~3. 
As we are going to outline, the finding of this proof was greatly
facilitated by the {\sl Mathematica} package HYP \cite{KratBF} developed by 
the first author. 
However, it must be noted that this proof yields
exclusively {\it formal\/} identities, i.e., identities for formal
power series in $x,y,z$; the identities are either meaningless or
wrong analytically. To remedy this fact,
in Section~4 we provide a second, direct 
approach to find 
$q$-analogues of the equality \eqref{GGR0}=\eqref{GGR1}=\eqref{GGR2}.
This approach is based on Ramanujan's $_1\psi_1$ summation formula.
As a result, we obtain $q$-analogues which are valid as analytic as
well as formal identities.
This proof yields, as a by-product, also a hypergeometric proof of the
identity of Gelfand, Graev and Retakh. 
However, it is a {\it genuine}
{\it basic} hypergeometric proof, i.e., it seems that the use of the base $q$
cannot be avoided in that proof.

\end{section}

\begin{section}{The proof}

We remark that the expressions \eqref{GGR1} and \eqref{GGR2} are
equal because of the transformation formula (cf.\ \cite[(1.7.1.3)]{SlatAC})
$$
{} _{2} F _{1} \!\left [ \begin{matrix} { a, b}\\ { c}\end{matrix} ; 
{\displaystyle
   z}\right ]  = {{{\left( 1 - z \right) }^{-a}}}
 {{{} _{2} F _{1} \!\left [ \begin{matrix} { a, c-b}\\ {
      c}\end{matrix} ; {\displaystyle -{\frac z {1 - z}}}\right ] }.
    }
$$
(In fact, the equality which is given explicitly in \cite{GeGRAA} is
the equality between \eqref{GGR0} and \eqref{GGR2}.)
Hence it would suffice to establish one of the equalities
\eqref{GGR0}=\eqref{GGR1} or \eqref{GGR0}=\eqref{GGR2}. 
However, since we eventually want to extend our proofs to the
$q$-case, we shall provide {\it direct\/} proofs for {\it both\/}
equalities.

We start with the proof of the equality \eqref{GGR0}=\eqref{GGR1}.
The strategy that we pursue is as follows: We compare coefficients of 
$x^Xy^Yz^Z$ in \eqref{GGR0} and
\eqref{GGR1}. The corresponding coefficient in \eqref{GGR0} is
the compact expression
\begin{equation} \label{eq:LS} 
{\frac{(\alpha)_X\,(\beta)_Y\,(1-\gamma)_Z\,(\gamma)_{X+Y-Z}} 
{X!\,Y!\,Z!\,(\alpha+\beta)_{X+Y-Z}}},
\end{equation}
while the corresponding coefficient in \eqref{GGR1} is the multiple sum
\begin{multline} \label{eq:RS} 
{\sum _{i,j,k,l,m,n} ^{}\kern-7pt}{\vphantom{\Big\vert}}'\kern7pt
 {\left( -1 \right) }^{i + j + k + l + m}
 \frac{       ({ \textstyle \beta}) _{n} \,
     ({ \textstyle \alpha + \beta - \gamma}) _{n} }
     {n!\,({ \textstyle \alpha + \beta}) _{n} }\\
\cdot
   {\binom{\beta - \gamma}  i}\,
     {\binom{-\beta - n}  j}\,
    {\binom{\alpha + \beta + n-1}  k}\,
     {\binom{\gamma -\alpha - \beta  - n}  l}\,
     {\binom n  m},
\end{multline}
where the sum ${\sum}'$ is subject to
\begin{equation} \label{eq:lin}
 i+l+n-m=X,\quad j+m=Y\quad \text{and}\quad k+l=Z.
\end{equation}
Because of the linear dependencies of the 6 parameters $i,j,k,l,m,n$,
the sum \eqref{eq:RS} can be written as a triple sum, once one
expresses all the parameters in terms of three fixed ones out of them.

The task is then to simplify this triple sum, by using known
hypergeometric summation and transformation formulas, interchange of
summations, and similar manipulations, until one arrives at the
compact expression \eqref{eq:LS}.
Clearly, in view of the many possibilities there are to choose 3 out
of 6 parameters, and the many possibilities to proceed afterwards,
this is a daunting task if one is to do this by hand. However, the
use of the computer will greatly facilitate the search for a 
feasible path. Indeed, by using the {\sl Mathematica} package HYP 
\cite{KratBF}, it did not take us more than three attempts to find
such a classical hypergeometric proof of the equality of \eqref{eq:LS}
and \eqref{eq:RS}. As a matter of fact, 
all the subsequent calculations were first
carried out on the computer by using HYP, and were then directly transformed
into \TeX-code using the {\sl Mathematica} command 
{\tt TeXForm}.\footnote{The corresponding {\sl Mathematica}
Notebook is available at:\newline {\tt
http://euler.univ-lyon1.fr/home/kratt/artikel/gelfand1.html}.} 

We are going to express all the parameters in terms of $i$, $j$ and
$k$. That is, we write
$$l= -k+Z,\quad m= -j+Y,\quad  n= -i-j+k+X+Y-Z,$$
and substitute this in \eqref{eq:RS}. Making the sum over $j$ the inner-most
sum, and writing it in hypergeometric notation, we obtain
\begin{multline*} \label{}
\sum _{i,k\ge0} ^{} \Bigg( {\left( -1 \right) }^{ Y + Z}\frac{
     ({ \textstyle \beta}) _{k + X + Y - Z-i} \,
     ({ \textstyle \alpha + \beta - \gamma}) _{
      k + X + Y - Z-i}  \,({ \textstyle  \gamma - \beta})
      _{i}}{i! \,
 k! \,Y! \,
     (k + X - Z-i)!  \,(Z-k)! }\\
\cdot
\frac {({ \textstyle 1 - \alpha - \beta +
      \gamma + i - X - Y}) _{Z-k}} 
{     ({ \textstyle \alpha + \beta}) _{X + Y - Z-i}}
     {} _{2} F _{1} \!\left [ \begin{matrix} {  1 - \alpha -
      \beta + i - X - Y + Z,-Y}\\ { 1 - \alpha - \beta
      + \gamma + i - X - Y}\end{matrix} ; {\displaystyle 1}\right ] \Bigg).
\end{multline*}
The $_2F_1$-series can be summed by means of the Chu--Vandermonde
summation 
(see \cite[(1.7.7), Appendix~(III.4)]{SlatAC}),
\begin{equation} \label{vand}
{} _{2} F _{1} \!\left [ \begin{matrix} { a, -N}\\ { c}\end{matrix} ;
   {\displaystyle 1}\right ]  = 
  {\frac{({ \textstyle c-a}) _{N} }{({ \textstyle c}) _{N} }},
\end{equation}
where $N$ is a nonnegative integer, because $Y$ is a nonnegative
integer. We substitute the result and now make the sum over $k$
the inner-most sum. If it is written in hypergeometric notation, 
we get
\begin{multline*} \label{}
\sum _{i\ge0} ^{} \Bigg(\frac{
     ({ \textstyle \beta}) _{X + Y - Z-i} \,
     ({ \textstyle \alpha + \beta - \gamma}) _{
      X -i} \,({ \textstyle \gamma - \beta})
      _{i} }{i! \,
   Y! \,Z! }\\
\cdot
\frac {
     ({ \textstyle \gamma - Z}) _{Y} } 
{(X - Z-i)! \,
      \,({ \textstyle \alpha + \beta})
      _{X + Y - Z-i}}
     {} _{2} F _{1} \!\left [ \begin{matrix} {  \beta - i + X +
      Y - Z,-Z}\\ { 1 - i + X - Z}\end{matrix} ; {\displaystyle 1}\right ]\Bigg) .
\end{multline*}
Again, the $_2F_1$-series can be summed by means of the Chu--Vandermonde
summation \eqref{vand}, since $Z$ is a nonnegative integer. 
Upon substituting the result, and writing the
remaining sum over $i$ in hypergeometric notation, we arrive at the 
expression
\begin{multline*} \label{}
  \frac{     ({ \textstyle \beta}) _{X + Y - Z} \,
     ({ \textstyle \alpha + \beta - \gamma}) _{X} \,
     ({ \textstyle 1 - \beta - Y}) _{Z} \,
     ({ \textstyle \gamma - Z}) _{Y} }{X! \,
   Y! \,Z! \,
     ({ \textstyle \alpha + \beta}) _{X + Y - Z} }\\
\times
{} _{3} F _{2} \!\left [ \begin{matrix} { 1 - \alpha -
      \beta - X - Y + Z, -\beta + \gamma, -X}\\ { 1
      - \beta - X - Y + Z, 1 - \alpha - \beta +
      \gamma - X}\end{matrix} ; {\displaystyle 1}\right ].
\end{multline*}
This time, the Pfaff--Saalsch\"utz summation formula
(see \cite[(2.3.1.3), Appendix~(III.2)]{SlatAC}),
\begin{equation} \label{pfaff}
{_3}F_2\!\left[\begin{matrix}a,b,-N\\
c,1+a+b-c-N\end{matrix};1\right]=
\frac{(c-a)_N\,(c-b)_N}{(c)_N\,(c-a-b)_N},
\end{equation}
where $N$ is a nonnegative integer, applies, because $X$ is a
nonnegative integer. Some simplification of
the result finally yields \eqref{eq:LS}, which finishes the
proof of the equality \eqref{GGR0}=\eqref{GGR1}.

Next we prove the equality \eqref{GGR0}=\eqref{GGR2}, following an
analogous approach. If we compare coefficients of $x^Xy^Yz^Z$
on both sides, we see
that we have to prove the equality of \eqref{eq:LS} with
\begin{multline} \label{eq:RS1} 
{\sum _{i,j,k,l,m,n} ^{}\kern-7pt}{\vphantom{\Big\vert}}'\kern7pt
{{\left( -1 \right) }^{i + j + k + l + m}}\,
  \frac {( \textstyle \beta)_n\,( \gamma) _{n}} {n!\,( \textstyle \alpha +
   \beta)_n} \\
\cdot
  {\binom {-\gamma - n} i}{\binom {-\beta - n} j}
  {\binom {\alpha + \beta + n-1} k}{\binom {\gamma-\alpha } l}
  {\binom n m},
\end{multline}
where the sum ${\sum}'$ is subject to
\begin{equation*} \label{eq:lin1}
 i+l+m=X,\quad j+n-m=Y\quad \text{and}\quad j+k+l=Z.
\end{equation*}
Here, we express all the parameters in terms of $i$, $k$ and
$l$. That is, we write
\begin{equation} \label{eq:lin2}
j = -k - l + Z,\quad m= -i - l + X,\quad  n= -i + k + X + Y - Z,
\end{equation}
and substitute this in \eqref{eq:RS1}. Making the sum over $k$ the inner-most
sum, and writing it in hypergeometric notation, we obtain
\begin{multline*} \label{}
\sum _{i,l\ge0} ^{} \Bigg( {{\left( -1 \right) }^{X+l+i}}{\frac {
     ({ \textstyle \beta}) _{X + Y-i - l} \,
     ({ \textstyle \alpha - \gamma}) _{l} \,
     ({ \textstyle \gamma}) _{X + Y - Z}} 
   {i!\,l!\,\left( X-i - l \right) !\,\left(  Y+l - Z \right) !\,
     \left( Z-l \right) !\,({ \textstyle \alpha + \beta}) _{ X + Y
      - Z-i} }}\\
\cdot
{} _{2} F _{1} \!\left [ \begin{matrix} {  \gamma + X + Y - Z,l - Z}\\ { 1 + l
      + Y - Z}\end{matrix} ; {\displaystyle 1}\right ]\Bigg).
\end{multline*}
The $_2F_1$-series can be summed by means of the Chu--Vandermonde
summation \eqref{vand}, because $Z-l$ is a nonnegative
integer. We substitute the result and now make the sum over $i$
the inner-most sum. If it is written in hypergeometric notation, 
we get
\begin{multline*} \label{}
\sum _{l\ge0} ^{} \Bigg(
{{\left( -1 \right) }^{X+l }}{\frac {
      ({ \textstyle \beta}) _{X + Y-l} \,
     ({ \textstyle \alpha - \gamma}) _{l} \,
     ({ \textstyle \gamma}) _{X + Y - Z} \,
     ({ \textstyle 1 - \gamma + l - X}) _{ Z-l} } 
   {l!\,\left( X-l \right) !\,Y!\,\left(  Z-l \right) !\,
     ({ \textstyle \alpha + \beta}) _{X + Y - Z} }}\\
\cdot
     {} _{2} F _{1} \!\left [ \begin{matrix} { 1 - \alpha - \beta - X - Y + Z, l
      - X}\\ { 1 - \beta + l - X - Y}\end{matrix} ; {\displaystyle 1}\right
]\Bigg).
\end{multline*}
Again, the $_2F_1$-series can be summed by means of the Chu--Vandermonde
summation \eqref{vand}, since $X-l$ is a nonnegative integer. 
Upon substituting the result, and writing the
remaining sum over $i$ in hypergeometric notation, we arrive at the 
expression
\begin{equation*} \label{}
{{\left( -1 \right) }^X}
{\frac {({ \textstyle \beta}) _{X + Y} \,
     ({ \textstyle \gamma}) _{X + Y - Z} \,
     ({ \textstyle 1 - \gamma - X}) _{Z} \,
     ({ \textstyle \alpha - Z}) _{X} }
   {X!\,Y! \,
     Z! \,({ \textstyle \alpha + \beta}) _{X + Y - Z}
      \,({ \textstyle 1 - \beta - X - Y}) _{X} }}
{} _{3} F _{2} \!\left [ \begin{matrix} { \alpha - 
      \gamma, -Z, -X}\\ { 1 - \gamma - X, \alpha - Z}\end{matrix} ;
      {\displaystyle 1}\right ].
\end{equation*}
As in the previous derivation, it is now 
the Pfaff--Saalsch\"utz summation formula \eqref{pfaff} which can be
applied, because $X$ is a
nonnegative integer. Some simplification of
the result finally yields \eqref{eq:LS}, which finishes the
proof of the equality \eqref{GGR0}=\eqref{GGR2}.

\end{section}


\begin{section}{Formal $q$-analogues}

Once a path of proof is found for a hypergeometric identity, it may
hint at a way to find a $q$-analogue. This principle constituted the
original motivation of Richard Askey to pose his problem, as we
already mentioned in the Introduction. Indeed, since the only
identities that we used in the derivation in Section~2 were the
Chu--Vandermonde summation and the Pfaff--Saalsch\"utz summation, 
it is not difficult to come up with a $q$-analogue. What one has to do
is to replace Pochhammer symbols by $q$-shifted factorials
$$(a;q)_n:=\begin{cases}
\prod _{i=0} ^{n-1}(1-aq^i), & n>0,\\ 
1, & n=0,\\
\prod_{i=0}^{-n-1}(1-aq^{n+i})^{-1}, & n<0,\end{cases}$$
binomials
by $q$-binomials, and finally to insert the ``right" powers of $q$. Since
there are two $q$-analogues of the Chu--Vandermonde summation formula,
there are in fact several possible $q$-analogues, the derivation of 
one of which we shall describe below. 
While we did not succeed to
find a $q$-analogue of the equality \eqref{GGR0}=\eqref{GGR1} which
looks equally elegant as the original identity, we did succeed for
the equality \eqref{GGR0}=\eqref{GGR2}.

We start by proving the following $q$-analogue of the
equality of \eqref{GGR0} and \eqref{GGR1}:
\begin{multline} \label{GGRq}
\sum _{j,k,m\ge0} ^{}
     q^{-\binom k2-\binom m2}
{\alpha}^{-m}\, {\beta}^{j - m}\,{\gamma}^{2 m-j-k}\,
\frac{
     ({\let \over / \def\frac#1#2{#1 / #2} \alpha}; q) _{j} \,
     ({\let \over / \def\frac#1#2{#1 / #2} \beta}; q) _{k} \,
    ({\let \over / \def\frac#1#2{#1 / #2}
      \frac{q}{\gamma}}; q) _{m}\,
   ({\let \over / \def\frac#1#2{#1 / #2} \gamma}; q) _{j + k -
      m}  }{
     ({\let \over / \def\frac#1#2{#1 / #2} q}; q) _{j} \,
     ({\let \over / \def\frac#1#2{#1 / #2} q}; q) _{k} \,
     ({\let \over / \def\frac#1#2{#1 / #2} q}; q) _{m} \,
({\let \over / \def\frac#1#2{#1 / #2}
      \alpha \beta}; q) _{j + k - m} }
x^jy^kz^m\\
= \sum _{i,j,k,l,m,n\ge0} ^{}
\Bigg(
 q^{mn
  - kl  - jm  - ln 
-\binom j2 - \binom k2 - \binom l2 - \binom m2}\,
{\alpha}^{-l}\,{\beta}^{i - l}\,{\gamma}^{-i - j + k + l}\,
\frac {({\let \over / \def\frac#1#2{#1 / #2} \beta}; q) _{n} \,
     ({\let \over / \def\frac#1#2{#1 / #2}
      \frac{\alpha \beta}{\gamma}}; q) _{n}} 
{     ({\let \over / \def\frac#1#2{#1 / #2}
      \alpha \beta}; q) _{n} \,
     ({\let \over / \def\frac#1#2{#1 / #2} q}; q) _{n}}\kern1.3cm\\
\frac{ 
     ({\let \over / \def\frac#1#2{#1 / #2}
      \frac{\gamma}{\beta}}; q) _{i} \,
     ({\let \over / \def\frac#1#2{#1 / #2} \beta q^ n};
      q) _{j} \,
 ({\let \over / \def\frac#1#2{#1 / #2} \frac{q^ {1 -
      n}}{\alpha \beta}}; q) _{k} \,
    ({\let \over / \def\frac#1#2{#1 / #2}
      \frac{\alpha \beta q^ n}{\gamma}}; q)
      _{l} \,
     ({\let \over / \def\frac#1#2{#1 / #2} q^ {-n}}; q) _{m}}{
     ({\let \over / \def\frac#1#2{#1 / #2} q}; q) _{i} \,
     ({\let \over / \def\frac#1#2{#1 / #2} q}; q) _{j} \,
     ({\let \over / \def\frac#1#2{#1 / #2} q}; q) _{k} \,
     ({\let \over / \def\frac#1#2{#1 / #2} q}; q) _{l} \,
     ({\let \over / \def\frac#1#2{#1 / #2} q}; q) _{m}  }\,
x^{i+l+n-m}y^{j+m}z^{k+l}\Bigg).
\end{multline}
The reader must be warned at this point that the only way to give this
identity a meaningful interpretation is as a formal power series in
$x,y,z$, and our proof will adopt this point of view. 
(Analytically, the series on the left-hand side of
\eqref{GGRq} diverges if $\vert q\vert<1$ because of the quadratic
powers of $q$. If $\vert q\vert>1$ the right-hand side diverges.)

For the proof of \eqref{GGRq} we proceed in complete analogy to
the previous section. It is needless to say that all the subsequent
calculations were again first carried out on the computer, this time
using the ``$q$-analogue" of HYP, HYPQ \cite{KratBF}, after which
they were transformed
into \TeX-code using the {\sl Mathematica} command 
{\tt TeXForm}.\footnote{The corresponding {\sl Mathematica}
Notebook is available at:\newline {\tt
http://euler.univ-lyon1.fr/home/kratt/artikel/gelfand1.html}.} 

As before, we compare coefficients of 
$x^Xy^Yz^Z$ on both sides of \eqref{GGRq}. We start with the
corresponding coefficient on the right-hand side. In analogy with
\eqref{eq:RS}, it is expressed 
in terms of a multiple sum over $i,j,k,l,m,n$ subject to
\eqref{eq:lin}.
We express all the 
parameters $i,j,k,l,m,n$ in terms of $i$, $j$ and $k$, make the
sum over $j$ the inner-most
sum, and write it in the standard basic hypergeometric notation 
$${}_{r+1}\phi_r\!\left[\begin{matrix} 
a_1,\dots,a_{r+1}\\ b_1,\dots,b_r\end{matrix}; q,
z\right]=\sum _{\ell=0} ^{\infty}\frac {\poq{a_1}{\ell}\cdots\poq{a_{r+1}}{\ell}}
{\poq{q}{\ell}\,\poq{b_1}{\ell}\cdots\poq{b_r}{\ell}}z^\ell\ .$$
Thus we obtain for the coefficient of $x^Xy^Yz^Z$ on 
the right-hand side of \eqref{GGRq} the expression
\begin{multline*} \label{}
\sum _{i,k\ge0} ^{} \Bigg(
{\left( -1 \right) }^{k + Y}\, q^{\binom {k+1}2+\binom {Z+1}2
 +Z( i - k - X - Y)}\,
{\alpha}^{-Z}\,{\beta}^{i - Z}\,
     {\gamma}^{Z-i}\,\\
\cdot
 \frac{  
     ({\let \over / \def\frac#1#2{#1 / #2} \beta}; q) _{k +
      X + Y - Z-i} \,({\let \over / \def\frac#1#2{#1 / #2}
      \frac{\alpha \beta}{\gamma}}; q) _{X +
      Y-i} \,({\let \over / \def\frac#1#2{#1 / #2}
      \frac{\gamma}{\beta}}; q) _{i} }{
     ({\let \over / \def\frac#1#2{#1 / #2} q}; q) _{i} \,
     ({\let \over / \def\frac#1#2{#1 / #2} q}; q) _{k} \,
     ({\let \over / \def\frac#1#2{#1 / #2} q}; q) _{Y} \,
     ({\let \over / \def\frac#1#2{#1 / #2} q}; q) _{k + X - Z-i} \,
     ({\let \over / \def\frac#1#2{#1 / #2} q}; q) _{Z-k} \,
 ({\let \over / \def\frac#1#2{#1 / #2}
      \alpha \beta}; q) _{X + Y - Z-i}  }\\
\cdot
{} _{2} \phi _{1} \! \left [
      \begin{matrix} \let \over / \def\frac#1#2{#1 / #2} \frac{q^
      {1 + i - X - Y + Z}}{\alpha \beta}, q^ {-Y}\\
      \let \over / \def\frac#1#2{#1 / #2} \frac{\gamma q^
      {1 + i - X - Y}}{\alpha \beta}\end{matrix} ;q,
      {\displaystyle q} \right ]
\Bigg).
\end{multline*}
The $_2\phi_1$-series can be summed by means of the following $q$-analogue of
the Chu--Vandermonde
summation (see \cite[(1.5.3); Appendix (II.6)]{GaRaAA}),
\begin{equation} \label{qvand} 
{}_2\phi _1\!\left [ \begin{matrix} \let\over/ a,{q^{-N}}\\ \let\over/  c\end{matrix} ;q,q\right ] =
  {\frac {{a^N} {(\let\over/ {c/ a};q)}_{N}} {{(\let\over/ c;q)}_{N}}},
\end{equation}
where $N$ is a nonnegative integer.
 We substitute the result and now make the sum over $k$
the inner-most sum. If it is written in basic hypergeometric notation, 
we get
\begin{multline*} \label{}
\sum _{i\ge0} ^{} \Bigg(
    q^{-\binom Y2+\binom {Z+1}2  + Z(i - X) }\,
{\alpha}^{-Z}\,{\beta}^{i - Z}\,{\gamma}^{-i - Y + Z}\,
\frac {  ({\let \over / \def\frac#1#2{#1 / #2} \beta}; q) _{ X +
      Y - Z-i} \,({\let \over / \def\frac#1#2{#1 / #2}
      \frac{\alpha \beta}{\gamma}}; q) _{X-i} }
 {   ({\let \over / \def\frac#1#2{#1 / #2} q}; q) _{i} \,
     ({\let \over / \def\frac#1#2{#1 / #2} q}; q) _{Y} \,
     ({\let \over / \def\frac#1#2{#1 / #2} q}; q) _{Z} }\\
\cdot
 \frac{ 
     ({\let \over / \def\frac#1#2{#1 / #2}
      \frac{\gamma}{\beta}}; q) _{i} \,
     ({\let \over / \def\frac#1#2{#1 / #2}
      \frac{\gamma}{q^ Z}}; q) _{Y} }{
     ({\let \over / \def\frac#1#2{#1 / #2} q}; q) _{X - Z-i} \,
   ({\let \over / \def\frac#1#2{#1 / #2}
      \alpha \beta}; q) _{X + Y - Z-i} }
     {} _{2} \phi _{1} \! \left [                \begin{matrix} \let \over /
      \def\frac#1#2{#1 / #2} \beta q^ {-i + X + Y - Z},
      q^ {-Z}\\ \let \over / \def\frac#1#2{#1 / #2} q^ {1
      - i + X - Z}\end{matrix} ;q, {\displaystyle q} \right ]\Bigg) .
\end{multline*}
Again, the $_2\phi_1$-series can be summed by means of the $q$-Chu--Vandermonde
summation \eqref{qvand}.
Upon substituting the result, and writing the
remaining sum over $i$ in basic hypergeometric notation, we arrive at the 
expression
\begin{multline*} \label{}
     q^{-\binom Y2-\binom Z2 + YZ }\,
{\alpha}^{-Z}\,{\gamma}^{-Y + Z}\,
\frac{
({\let \over / \def\frac#1#2{#1 / #2} \beta}; q)
      _{X + Y - Z} \,({\let \over / \def\frac#1#2{#1 / #2}
      \frac{\alpha \beta}{\gamma}}; q) _{X} \,
     ({\let \over / \def\frac#1#2{#1 / #2} \frac{q^ {1 -
      Y}}{\beta}}; q) _{Z} \,
     ({\let \over / \def\frac#1#2{#1 / #2}
      \frac{\gamma}{q^ Z}}; q) _{Y} }{
     ({\let \over / \def\frac#1#2{#1 / #2} q}; q) _{X} \,
     ({\let \over / \def\frac#1#2{#1 / #2} q}; q) _{Y} \,
     ({\let \over / \def\frac#1#2{#1 / #2} q}; q) _{Z} \,
     ({\let \over / \def\frac#1#2{#1 / #2}
      \alpha \beta}; q) _{X + Y - Z} }\\
\times
      {} _{3} \phi _{2} \! \left [           
     \begin{matrix} \let \over /      \def\frac#1#2{#1 / #2}
         \frac{q^ {1 - X - Y + Z}}{\alpha \beta},
     \frac{\gamma}{\beta},
   q^ {-X}\\ \let \over / \def\frac#1#2{#1 / #2}
      \frac{q^ {1 - X - Y + Z}}{\beta},      \frac{\gamma q^ {1 -
      X}}{\alpha \beta}\end{matrix} ;q, {\displaystyle q}
      \right ].
\end{multline*}
The $_3\phi_2$-series can be summed by means of the $q$-analogue of
the Pfaff--Saalsch\"utz summation formula
(see \cite[(1.7.2); Appendix (II.12)]{GaRaAA}),
\begin{equation} \label{qpfaff}
{}_3\phi _2\!\left [ \begin{matrix} \let\over/ a,b,{q^{-N}}\\ \let\over/  c,{{a b {q^{1 - N}}}\over
   c}\end{matrix} ;q,q\right ] = {\frac{{(\let\over/ {c/ a};q)}_{N} {(\let\over/ {c/b};q)}_{N}}
    {{(\let\over/ c;q)}_{N} {(\let\over/ {c/ {a b}};q)}_{N}}},
\end{equation}
where $N$ is a nonnegative integer. Some simplification then yields
the expression 
$$     q^{-\binom Y2-\binom Z2}
{\alpha}^{-Z}\, {\beta}^{X - Z}\,{\gamma}^{2Z-X-Y}\,
\frac{
     ({\let \over / \def\frac#1#2{#1 / #2} \alpha}; q) _{X} \,
     ({\let \over / \def\frac#1#2{#1 / #2} \beta}; q) _{Y} \,
    ({\let \over / \def\frac#1#2{#1 / #2}
      \frac{q}{\gamma}}; q) _{Z}\,
   ({\let \over / \def\frac#1#2{#1 / #2} \gamma}; q) _{X + Y -
      Z}  }{
     ({\let \over / \def\frac#1#2{#1 / #2} q}; q) _{X} \,
     ({\let \over / \def\frac#1#2{#1 / #2} q}; q) _{Y} \,
     ({\let \over / \def\frac#1#2{#1 / #2} q}; q) _{Z} \,
({\let \over / \def\frac#1#2{#1 / #2}
      \alpha \beta}; q) _{X + Y - Z} },
$$
which is exactly the coefficient of $x^Xy^Yz^Z$ on the left-hand side
of \eqref{GGRq}, thus establishing our claim.

In a completely analogous manner, we can prove the following
$q$-analogue of the 
equality of \eqref{GGR0} and \eqref{GGR2}:
\begin{multline} \label{GGRq1}
\sum _{j,k,m\ge0} ^{}
{\frac {({\let \over / \def\frac#1#2{#1 / #2} \alpha}; q) _{j} \,
     ({\let \over / \def\frac#1#2{#1 / #2} \beta}; q) _{k} \,
     ({\let \over / \def\frac#1#2{#1 / #2} {q\over {\gamma}}}; q) _{m} \,
     ({\let \over / \def\frac#1#2{#1 / #2} \gamma}; q) _{j + k - m} }
   {({\let \over / \def\frac#1#2{#1 / #2} q}; q) _{j} \,
     ({\let \over / \def\frac#1#2{#1 / #2} q}; q) _{k} \,
     ({\let \over / \def\frac#1#2{#1 / #2} q}; q) _{m} \,
     ({\let \over / \def\frac#1#2{#1 / #2} \alpha \beta};
      q) _{j + k - m} }}
x^jy^kz^m\\
= \sum _{i,j,k,l,m,n\ge0} ^{}
\Bigg(
{q^{ k n- i n +m }}\,{{\alpha}^k}\,{{\beta}^{k -i-m}}{{\gamma}^{-k}}
\frac {
     ({\let \over / \def\frac#1#2{#1 / #2} \beta}; q) _{n} \,
({\let \over / \def\frac#1#2{#1 / #2} \gamma}; q) _{n} } 
{
({\let \over / \def\frac#1#2{#1 / #2} \alpha \beta};
      q) _{n} \,
     ({\let \over / \def\frac#1#2{#1 / #2} q}; q) _{n} }\kern5.5cm\\
\cdot
\frac {
     ({\let \over / \def\frac#1#2{#1 / #2} \gamma {q^n}}; q) _{i} \,
     ({\let \over / \def\frac#1#2{#1 / #2} \beta {q^n}}; q) _{j}\,
     ({\let \over / \def\frac#1#2{#1 / #2} {{{q^{1 - n}}}\over {\alpha {
      \beta}}}}; q) _{k} \,
    ({\let \over / \def\frac#1#2{#1 / #2} {{\alpha}\over {\gamma}}}; q)
      _{l} \, ({\let \over / \def\frac#1#2{#1 / #2} {q^{-n}}}; q)
      _{m} } 
   {({\let \over / \def\frac#1#2{#1 / #2} q}; q) _{i} \,
     ({\let \over / \def\frac#1#2{#1 / #2} q}; q) _{j} \,
     ({\let \over / \def\frac#1#2{#1 / #2} q}; q) _{k} \,
     ({\let \over / \def\frac#1#2{#1 / #2} q}; q) _{l} \,
     ({\let \over / \def\frac#1#2{#1 / #2} q}; q) _{m} }
x^{i+l+m}y^{j+n-m}z^{j+k+l}\Bigg).
\end{multline}
That is, as in the proof of the equality \eqref{GGR0}=\eqref{GGR2} in
the previous section, we compare coefficients of $x^Xy^Yz^Z$ on both
sides, then start with the the coefficient of the right-hand side,
substitute the relations \eqref{eq:lin2}, evaluate the sums over
$k$ and $i$ (in that order) by means of the $q$-analogue \eqref{qvand}
of the Chu--Vandermonde
summation, and finally evaluate the sum over $l$ by
means of the the $q$-analogue \eqref{qpfaff} of the
Pfaff--Saalsch\"utz summation. 
However, again, this is an identity which makes sense only
as a formal power series in $x,y,z$. 
(Analytically, the right-hand side never converges.)

Remarkably, in this case, the right-hand side of \eqref{GGRq1} can be
simplified. To be precise, the sums over $i,j,k,l,m$ can be evaluated
with the help of the $q$-binomial theorem (see \cite[(1.3.2); Appendix
(II.3)]{GaRaAA})
\begin{equation} \label{qbin}
\sum _{\ell=0} ^{\infty}\frac {(a;q)_\ell} {(q;q)_\ell}z^\ell
=  \frac {{( \let\over/ a z;q)_\infty}}
{{( \let\over/  z;q)_\infty}}.
\end{equation}
If the remaining sum over $n$ is written in basic hypergeometric
notation, the result is the compact identity
\begin{multline} \label{GGRq2}
\sum _{j,k,m\ge0} ^{}
{\frac {({\let \over / \def\frac#1#2{#1 / #2} \alpha}; q) _{j} \,
     ({\let \over / \def\frac#1#2{#1 / #2} \beta}; q) _{k} \,
     ({\let \over / \def\frac#1#2{#1 / #2} {q\over {\gamma}}}; q) _{m} \,
     ({\let \over / \def\frac#1#2{#1 / #2} \gamma}; q) _{j + k - m} }
   {({\let \over / \def\frac#1#2{#1 / #2} q}; q) _{j} \,
     ({\let \over / \def\frac#1#2{#1 / #2} q}; q) _{k} \,
     ({\let \over / \def\frac#1#2{#1 / #2} q}; q) _{m} \,
     ({\let \over / \def\frac#1#2{#1 / #2} \alpha \beta};
      q) _{j + k - m} }}
x^jy^kz^m\\
= 
\frac {({\let \over / \def\frac#1#2{#1 / #2} {{\gamma x}\over \beta}}; q)
      _{\infty} } 
{({\let \over / \def\frac#1#2{#1 / #2} {x\over \beta}}; q) _{\infty} }
\frac {({\let \over / \def\frac#1#2{#1 / #2} {{q z}\over {{
      \gamma}}}}; q) _{\infty} } 
{({\let \over / \def\frac#1#2{#1 / #2} {{\alpha \beta z}\over {{
      \gamma}}}}; q) _{\infty} }
\frac {({\let \over / \def\frac#1#2{#1 / #2} {{{
      \alpha} x z}\over {\gamma }}}; q) _{\infty} } 
{({\let \over / \def\frac#1#2{#1 / #2} {{{
      } x z}}}; q) _{\infty} }
{\frac {  ({\let \over / \def\frac#1#2{#1 / #2} \beta y z}; q) _{\infty} }
     {  ({\let \over / \def\frac#1#2{#1 / #2} y z}; q) _{\infty} }}\\
\times
{} _{4} \phi _{3} \! \left [             \begin{matrix} \let \over /
      \def\frac#1#2{#1 / #2} \beta, {{\alpha \beta z}\over {{
      \gamma}}}, {{\beta y}\over x}, \gamma\\ \let \over / \def\frac#1#2{#1 / #2}
      \alpha \beta, \beta y z, {{{q\beta}}\over x}\end{matrix} ;q,
      {\displaystyle q} \right ] ,
\end{multline}
which is a much more elegant $q$-analogue of the equality 
\eqref{GGR0}=\eqref{GGR2} than the identity \eqref{GGRq} is a
$q$-analogue of the equality \eqref{GGR0}=\eqref{GGR1}. (We remark that
the $_4\phi_3$-series on the right-hand side is balanced, i.e., the
product of the lower parameters is equal to $q$ times the product of
the upper parameters.) 
We repeat that this is an identity for {\it formal power series in
$x,y,z$}. It is in fact wrong (!) if one would interpret the left- and
right-hand sides as analytic series. Although both sides make perfect
sense analytically, on the right-hand side there is a term missing as
we are going to show in the next section (cf.\ \eqref{GGRq5}). The
problem, when the derivation of \eqref{GGRq2} is regarded in the
analytic sense, arises when in
\eqref{GGRq1} we use the $q$-binomial theorem \eqref{qbin} to evaluate
the sum over $i$. For, we would have to apply \eqref{qbin} with
$z=x/q^{n}$. If $\vert q\vert<1$, then
the left-hand side of \eqref{qbin}
converges only if $\vert z\vert<1$, but we cannot have $\vert
x/q^{n}\vert<1$ for arbitrarily large $n$ (unless $x=0$). There are
similar problems if $\vert q\vert>1$.

In the next section, we
shall not only derive a $q$-analogue of \eqref{GGR0}=\eqref{GGR2}
which does not have these problems, i.e., which 
is valid in the formal {\it as well as} in the analytic sense, but
as well a more elegant $q$-analogue of \eqref{GGR0}=\eqref{GGR1},
which is also valid in both senses.
\end{section}

\begin{section}{Analytic $q$-analogues}

In this section we outline a different method to obtain 
$q$-analogues of the equality \eqref{GGR0}=\eqref{GGR1}=\eqref{GGR2}. It is
based on Ramanujan's ${}_1\psi_1$ summation (see \cite[Ex.~1.6(ii); Appendix
(II.5)]{GaRaAA})
\begin{equation}\label{1psi1}
\sum_{n=-\infty}^\infty\frac{(a;q)_n}{(b;q)_n}z^n
=\frac{(q;q)_\infty\,(b/a;q)_\infty\,(az;q)_\infty\,(q/az;q)_\infty}
{(b;q)_\infty\,(q/a;q)_\infty\,(z;q)_\infty\,(b/az;q)_\infty}, 
\qquad |b/a|<|z|<1, 
\end{equation}
or equivalently,
$$\frac{(a;q)_n}{(b;q)_n}=\frac 1{2\pi i}\int
\frac{(q;q)_\infty\,(b/a;q)_\infty\,(az;q)_\infty\,(q/az;q)_\infty}
{(b;q)_\infty\,(q/a;q)_\infty\,(z;q)_\infty\,(b/az;q)_\infty}\,
z^{-n-1}\,dz,
$$
where the integration is over a contour encircling the origin
counter-clockwise inside the annulus of convergence. 
Thus,
\begin{align*}
\sum_{j,k,m\ge0}&\frac{(a;q)_j\,(b;q)_k\,(c;q)_m\,(e;q)_{j+k-m}}
{(q;q)_j\,(q;q)_k\,(q;q)_m
\,(f;q)_{j+k-m}}\,x^jy^kz^m\\
&=
\sum_{j,k,m\ge0}\frac{(a;q)_j\,(b;q)_k\,(c;q)_m}{(q;q)_j\,(q;q)_k\,(q;q)_m
}\,x^jy^kz^m\\
&\kern3cm
\cdot
\frac 1{2\pi i}\int
\frac{(q;q)_\infty\,(f/e;q)_\infty\,(et;q)_\infty\,(q/et;q)_\infty}
{(f;q)_\infty\,(q/e;q)_\infty\,(t;q)_\infty\,(f/et;q)_{\infty}}
\,t^{m-j-k-1}\,dt\\
&=\frac 1{2\pi i}\int
\frac{(q;q)_\infty\,(f/e;q)_\infty\,(et;q)_\infty\,
(q/et;q)_\infty\,(ax/t;q)_\infty\,(by/t;q)_\infty\,(czt;q)_\infty}
{(f;q)_\infty\,(q/e;q)_\infty\,(t;q)_\infty\,(f/et;q)_\infty\,
(x/t;q)_\infty\,(y/t;q)_\infty\,(zt;q)_{\infty}}
\,\frac{dt}{t},
\end{align*}
where we used the $q$-binomial theorem \eqref{qbin} to evaluate the
sums over $j$, $k$, and $m$. 
This requires the convergence conditions
$\max(|x|,|y|,|f/e|)<\min(1,1/|z|)$.
If we in addition assume $|ce|<1$,
 the value of the integral is given by \cite[(4.10.9)]{GaRaAA} as
\begin{multline*}
\frac{(ax;q)_\infty\,(by;q)_\infty\,(cz;q)_\infty\,(e;q)_\infty}
{(x;q)_\infty\,(y;q)_\infty\,(z;q)_\infty\,(f;q)_\infty}\,
{}_{4}\phi_3\left[\begin{matrix} {x,y,q/cz,f/e}\\
{ax,by,q/z}\end{matrix};q,ce\right]\\
+\frac{(axz;q)_\infty\,(byz;q)_\infty\,(c;q)_\infty\,
(e/z;q)_\infty\,(qz/e;q)_\infty\,(f/e;q)_\infty}
{(xz;q)_\infty\,(yz;q)_\infty\,(1/z;q)_\infty\,(f;q)_\infty\,
(q/e;q)_\infty\,(fz/e;q)_\infty}\\
\times
{}_{4}\phi_3\left[\begin{matrix}{xz,yz,q/c,fz/e}\\
{axz,byz,qz}\end{matrix};q,ce\right].
\end{multline*}

Until now, nothing ``special'' has been used. Indeed, using
\eqref{1psi1} instead of the $q$-binomial theorem, we get the integral
representation 
\begin{multline*}\sum_{k_1,\dots,k_n=-\infty}^{\infty}\frac{(a;q)_{k_1+\dots+k_n}}
{(c;q)_{k_1+\dots +k_n}}\prod_{j=1}^n
\frac{(b_j;q)_{k_j}}{(d_j;q)_{k_j}}\,x_j^{k_j}\\
=
\frac 1{2\pi i}\int
\frac{(q;q)_\infty\,(c/a;q)_\infty\,(at;q)_\infty\,
(q/at;q)_\infty}
{(c;q)_\infty\,(q/a;q)_\infty\,(t;q)_\infty\,(c/at;q)_\infty}\\
\times
\prod_{j=1}^n\frac{(q;q)_\infty\,(d_j/b_j;q)_\infty\,(b_jx_j/t;q)_\infty\,
(qt/b_jx_j;q)_\infty}{(d_j;q)_\infty\,(q/b_j;q)_\infty\,(x_j/t;q)_\infty\,
(d_jt/b_jx_j;q)_\infty}\,\frac{dt}{t}, \end{multline*}
where 
 $\max(c/a,x_1,\dots,x_n)<|t|<\min(1,x_1b_1/d_1,\dots,x_nb_n/d_n)$
on the contour of integration. This type of integral may be expressed
as a finite sum of basic hypergeometric series, cf.\ 
\cite[(4.10.8), (4.10.9)]{GaRaAA}.
The case considered above is the special case $n=3$, $d_1=d_2=q$,
$b_3=1$. The case $d_1=\dots=d_n=q$ gives Andrews's identity 
\cite{andrews}
\begin{multline*}
\sum_{k_1,\dots,k_n=0}^{\infty}\frac{(a;q)_{k_1+\dots+k_n}}
{(c;q)_{k_1+\dots +k_n}}\prod_{j=1}^n
\frac{(b_j;q)_{k_j}}{(q;q)_{k_j}}\,x_j^{k_j}\\
=\frac{(c;q)_\infty}{(d;q)_\infty}
\prod_{j=1}^n\frac{(a_jx_j;q)_\infty}{(x_j;q)_\infty}\,
{}_{n+1}\phi_n\left[\begin{matrix} {d/c,x_1,\dots,x_n}\\
{ax_1,\dots,ax_n}\end{matrix};q,c\right].
\end{multline*}

Returning to the special case under consideration,
 we choose $ab=f$ and $ce=q$. This is the condition that the
${}_4\phi_3$-series 
are balanced (i.e., that the
product of the lower parameters is equal to $q$ times the product of
the upper parameters). They may then be 
combined to a very-well-poised ${}_8\phi_7$-series using 
\cite[(2.10.10), Appendix (III.36)]{GaRaAA}. 
After replacing $a$ by $\alpha$, $b$ by $\beta$, and $e$ by $\gamma$,
the conclusion is that
\begin{multline} \label{eq:8phi7}
\sum_{j,k,m\ge0}\frac{(\alpha ;q)_j\,(\beta ;q)_k\,(q/\gamma ;q)_m\,(\gamma ;q)_{j+k-m}}
{(q;q)_j\,(q;q)_k\,(q;q)_m
\,(\alpha \beta ;q)_{j+k-m}}\,x^jy^kz^m\\
=
\frac {(\alpha x;q)_\infty} {(x;q)_\infty}
\frac {(\beta y;q)_\infty} {(y;q)_\infty}
\frac {(\gamma y;q)_\infty} {(\alpha \beta y;q)_\infty}
\frac {(qz/\gamma ;q)_\infty} {(\alpha \beta z/\gamma ;q)_\infty}
\frac{(\alpha \beta yz/\gamma ;q)_\infty}
{(yz;q)_\infty}\kern4cm\\
\times
{} _{8} \phi _{7} \! \left [             \begin{matrix} \let \over /
      \def\frac#1#2{#1 / #2} {{\alpha \beta y}\over q}, {\sqrt{{
      \alpha} \beta q y}}, -{\sqrt{\alpha \beta q y}}, \alpha, y,
      {{\gamma}\over z}, {{\alpha \beta}\over {\gamma}}, {{{
      \beta} y}\over x}\\ \let \over / \def\frac#1#2{#1 / #2} 
{{{\sqrt{\alpha {
      \beta} y}}}\over {{\sqrt{q}}}}, -{{{\sqrt{\alpha {
      \beta} y}}}\over {{\sqrt{q}}}}, \beta y, \alpha \beta, {{{
      \alpha} \beta y z}\over {\gamma}}, \gamma y, {
      \alpha} x\end{matrix} ;q, {\displaystyle x z} \right ] 
.
\end{multline}
Finally, we apply Bailey's
very-well-poised $_8\phi_7$
transformation (see \cite[(2.10.1); Appendix (III.23)]{GaRaAA})
\begin{multline} \label{Bailey}
{} _{8} \phi _{7} \! \left [             \begin{matrix} \let \over / \def\frac#1#2{#1
   / #2} a, {\sqrt{a}} q, - {\sqrt{a}} q  , b, c, d, e, f\\
   \let \over / \def\frac#1#2{#1 / #2} {\sqrt{a}}, -{\sqrt{a}}, {{a q}\over
   b}, {{a q}\over c}, {{a q}\over d}, {{a q}\over e}, {{a q}\over
   f}\end{matrix} ;q, {\displaystyle {\frac {{a^2} {q^2}} {b c d e f}}}
   \right ] \\
= {\frac { ({\let \over / \def\frac#1#2{#1 / #2} a q}; q) _{\infty} \,
      ({\let \over / \def\frac#1#2{#1 / #2} {{a q}\over {e f}}}; q)
       _{\infty} \,({\let \over / \def\frac#1#2{#1 / #2} {{{a^2} {q^2}}\over
       {b c d e}}}; q) _{\infty} \,
      ({\let \over / \def\frac#1#2{#1 / #2} {{{a^2} {q^2}}\over
       {b c d f}}}; q) _{\infty} }
    {({\let \over / \def\frac#1#2{#1 / #2} {{a q}\over e}}; q) _{\infty} \,
      ({\let \over / \def\frac#1#2{#1 / #2} {{a q}\over f}}; q) _{\infty} \,
      ({\let \over / \def\frac#1#2{#1 / #2} {{{a^2} {q^2}}\over {b c d}}};
       q) _{\infty} \,({\let \over / \def\frac#1#2{#1 / #2}
       {{{a^2} {q^2}}\over {b c d e f}}}; q) _{\infty} }}\kern4cm\\
\times
{} _{8} \phi _{7} \! \left [             \begin{matrix} \let \over /
       \def\frac#1#2{#1 / #2} {{{a^2} q}\over {b c d}}, {{a {q^{{3\over
       2}}}}\over {{\sqrt{b c d}}}}, -{{a {q^{{3\over 2}}}}\over
       {{\sqrt{b c d}}}}, {{a q}\over {c d}}, {{a q}\over {b d}},
       {{a q}\over {b c}}, e, f\\ \let \over / \def\frac#1#2{#1 / #2}
       {{a {\sqrt{q}}}\over {{\sqrt{b c d}}}}, -{{a {\sqrt{q}}}\over
       {{\sqrt{b c d}}}}, {{a q}\over b}, {{a q}\over c}, {{a q}\over d},
       {{{a^2} {q^2}}\over {b c d e}}, {{{a^2} {q^2}}\over
       {b c d f}}\end{matrix} ;q, {\displaystyle {\frac {a q} {e f}}} \right
       ].
\end{multline}
As a result, we obtain the identity
\begin{multline} \label{GGRq3}
\sum_{j,k,m\ge0}\frac{(\alpha ;q)_j\,(\beta ;q)_k\,(q/\gamma ;q)_m\,(\gamma ;q)_{j+k-m}}
{(q;q)_j\,(q;q)_k\,(q;q)_m
\,(\alpha \beta ;q)_{j+k-m}}\,x^jy^kz^m\\
=
\frac {({\let \over / \def\frac#1#2{#1 / #2} {{\gamma x}\over
      {\beta}}}; q) _{\infty} } 
  {({\let \over / \def\frac#1#2{#1 / #2} x}; q) _{\infty} }
\frac {({\let \over / \def\frac#1#2{#1 / #2} \beta y}; q) _{\infty} } 
  {({\let \over / \def\frac#1#2{#1 / #2} y}; q) _{\infty} }
\frac {({\let \over / \def\frac#1#2{#1 / #2} {{q z}\over {\gamma}}}; q)
      _{\infty}} 
  {({\let \over / \def\frac#1#2{#1 / #2} {{\alpha \beta z}\over {{
      \gamma}}}}; q) _{\infty}}
\frac {     ({\let \over / \def\frac#1#2{#1 / #2} {{\alpha {
      \beta} x z}\over {\gamma}}}; q) _{\infty} } 
   {     ({\let \over / \def\frac#1#2{#1 / #2} x z}; q)
      _{\infty} }
{\frac { 
     ({\let \over / \def\frac#1#2{#1 / #2} \beta y z}; q) _{\infty} \,
     ({\let \over / \def\frac#1#2{#1 / #2} {{\alpha \beta y z}\over
      {\gamma}}}; q) _{\infty} }
   { ({\let \over / \def\frac#1#2{#1 / #2} y z}; q) _{\infty} \,
     ({\let \over / \def\frac#1#2{#1 / #2} {{\alpha {{{
      \beta}}^2} y z}\over {\gamma}}}; q) _{\infty} }}\kern2cm\\
\times
{} _{8} \phi _{7} \! \left [             \begin{matrix} \let \over /
      \def\frac#1#2{#1 / #2} {{\alpha {{\beta}^2} y z}\over {{
      \gamma} q}}, {{\beta {\sqrt{\alpha q y z}}}\over {{\sqrt{{
      \gamma}}}}}, -{{\beta {\sqrt{\alpha q y z}}}\over {{\sqrt{{
      \gamma}}}}}, {{\alpha \beta z}\over {\gamma}}, {{{
      \beta} y z}\over {\gamma}}, \beta, {{\alpha \beta}\over {{
      \gamma}}}, {{\beta y}\over x}\\ \let \over / \def\frac#1#2{#1 / #2}
      {{\beta {\sqrt{\alpha y z}}}\over {{\sqrt{\gamma q}}}},
      -{{\beta {\sqrt{\alpha y z}}}\over {{\sqrt{\gamma q}}}}, {
      \beta} y, \alpha \beta, {{\alpha \beta y z}\over {{
      \gamma}}}, \beta y z, {{\alpha \beta x z}\over {{
      \gamma}}}\end{matrix} ;q, {\displaystyle {\frac {\gamma x} {\beta}}}
      \right ]
,
\end{multline}
valid if $\max(|x|,|y|,|\alpha\beta/\gamma|)<\min(1,1/|z|)$ and
$|\gamma x/\beta|<1$, 
which is a perfect $q$-analogue of the equality
\eqref{GGR0}=\eqref{GGR1}. 
It is valid both analytically and as a
formal power series in $x,y,z$.

On the other hand, if we apply the transformation formula
\eqref{Bailey} to the $_8\phi_7$-series in \eqref{eq:8phi7} where the
lower and upper parameters are in the order
$$
{} _{8} \phi _{7} \! \left [             \begin{matrix} \let \over / \def\frac#1#2{#1
   / #2} {{\alpha \beta y}\over q}, {\sqrt{\alpha \beta q y}},
   -{\sqrt{\alpha \beta q y}}, \alpha, {{\alpha \beta}\over
   {\gamma}}, y, {{\beta y}\over x}, {{\gamma}\over z}\\ \let \over /
   \def\frac#1#2{#1 / #2} {{{\sqrt{\alpha \beta y}}}\over
   {{\sqrt{q}}}}, -{{{\sqrt{\alpha \beta y}}}\over {{\sqrt{q}}}}, {
   \beta} y, \gamma y, \alpha \beta, \alpha x, {{\alpha {
   \beta} y z}\over {\gamma}}\end{matrix} ;q, {\displaystyle x z} \right
],
$$
then the result is the identity
\begin{multline} \label{GGRq4}
\sum_{j,k,m\ge0}\frac{(\alpha ;q)_j\,(\beta ;q)_k\,(q/\gamma ;q)_m\,(\gamma ;q)_{j+k-m}}
{(q;q)_j\,(q;q)_k\,(q;q)_m
\,(\alpha \beta ;q)_{j+k-m}}\,x^jy^kz^m\\
=
\frac {( \gamma x;q)_\infty} {( x;q)_\infty}
\frac {\let\over/ ( {{q z}\over {\gamma}} ;q) _\infty} 
{\let\over/ ( {{{
   \alpha} \beta z}\over {\gamma}} ;q) _\infty}
\frac {\let\over/ ({{ \alpha x z}\over {\gamma}};q)_\infty} 
 {( x z;q)_\infty}
\frac {( \beta y z;q)_\infty} {( y z;q)_\infty}
 \frac {( \beta y;q)_\infty\,( {
   \gamma} y;q)_\infty} 
{\let\over/ (\beta \gamma y;q)_\infty\,( y;q)_\infty} \\
\times
{} _{8} \phi _{7} \! \left [ \begin{matrix} \let \over / \def\frac#1#2{#1
   / #2} {{\beta \gamma y}\over q}, {\sqrt{\beta \gamma q y}},
   -{\sqrt{\beta \gamma q y}}, \gamma, \beta, {{\gamma y}\over
   {\alpha}}, {{\beta y}\over x}, {{\gamma}\over z}\\ \let \over /
   \def\frac#1#2{#1 / #2} {{{\sqrt{\beta \gamma y}}}\over
   {{\sqrt{q}}}}, -{{{\sqrt{\beta \gamma y}}}\over {{\sqrt{q}}}}, {
   \beta} y, \gamma y, \alpha \beta, \gamma x, {
   \beta} y z\end{matrix} ;q, {\displaystyle {\frac {\alpha x z} {\gamma}}}
   \right ]
,
\end{multline}
which is another $q$-analogue of the equality
\eqref{GGR0}=\eqref{GGR2}. 
Again, it is not only valid analytically, but also as a
formal power series in $x,y,z$.

We remark that, by another application of
\cite[(III.36)]{GaRaAA}, the left-hand side of \eqref{GGRq4}
can alternatively be written
\begin{multline}\label{GGRq5}
\frac{(\gamma
  x/\beta;q)_\infty\,
(qz/\gamma;q)_\infty\,(\alpha xz/\gamma;q)_\infty\,(\beta
  yz;q)_\infty}{(x/\beta;q)_\infty\,
(\alpha\beta z/\gamma;q)_\infty\,(xz;q)_\infty\,
(yz;q)_\infty}\,{}_{4}\phi_3\left[
\begin{matrix}{\beta,\gamma,\alpha\beta
  z/\gamma,\beta y/x}\\
{\alpha\beta,\beta yz,q\beta /x}\end{matrix};q,q\right]\\
+\frac{(qz/\gamma;q)_\infty\,(\beta;q)_\infty\,
(\gamma;q)_\infty\,(\beta y/x;q)_\infty\,(\alpha
  x;q)_\infty\,(xyz;q)_\infty}
{(x;q)_\infty\,(y;q)_\infty\,(xz;q)_\infty\,(yz;q)_\infty\,
(\alpha\beta;q)_\infty\,(\beta/x;q)_\infty}\\
\times{}_{4}\phi_3
\left[\begin{matrix}{x,y,\gamma x/\beta,
 \alpha xz/\gamma}\\
{\alpha x,xyz,qx/\beta }\end{matrix};q,q\right].
\end{multline}
Note that the first term is precisely the (analytically false)
expression given in \eqref{GGRq2}.

An interesting aspect of this proof is that it provides in particular
a proof of the equality \eqref{GGR0}=\eqref{GGR1}=\eqref{GGR2} 
(by doing the replacements
$\alpha\to q^\alpha$, $\beta\to q^\beta$, $\gamma\to q^\gamma$, and
then performing the limit $q\to1$). However, because of the use of Ramanujan's
${}_1\psi_1$ summation, it seems impossible to make it into a proof ``with
$q=1$," i.e., into a proof that uses only identities for {\it
ordinary} hypergeometric series. On the other hand, we have given
such a proof in Section~2.

\end{section}

\end{document}